\documentclass[a4paper,dvipsnames]{article}

\usepackage{amsmath}
\usepackage{amssymb}
\usepackage{amsthm}
\usepackage{bm}
\usepackage{cite}
\usepackage[colorlinks=true,allcolors=purple]{hyperref}
\usepackage{cleveref}
	\crefformat{equation}{#2(#1)#3}
\usepackage{dsfont}
\usepackage{enumitem}
\usepackage[OT2,OT1]{fontenc}
\usepackage[cal=cm,scr=boondoxo]{mathalfa}
\usepackage{pifont}
\usepackage{stmaryrd}
\usepackage{textcomp}
\usepackage{tikz}
	\usetikzlibrary{arrows}
	\usetikzlibrary{patterns}
\usepackage{tikz-cd}
\usepackage[textwidth=3cm,colorinlistoftodos]{todonotes}
\usepackage{wasysym}
\usepackage{xfrac}

\numberwithin{equation}{section}			
\swapnumbers						

\newcommand\cyr{
\renewcommand\rmdefault{wncyr}%
\renewcommand\sfdefault{wncyss}%
\renewcommand\encodingdefault{OT2}%
\normalfont
\selectfont}
\DeclareTextFontCommand{\textcyr}{\cyr}

\newcounter{Enum}				
\newenvironment{Enumerate}{\begin{enumerate}[label={\rm({\roman*})}]}{\end{enumerate}}
\newcommand{\Enumref}[1]{{\setcounter{Enum}{#1}{\rm(\roman{Enum})}}}

\newcommand{\descriptionlabelsave}{}		
\newenvironment{Itemize}{%
	\renewcommand{\descriptionlabelsave}{\descriptionlabel}\renewcommand{\descriptionlabel}{$\triangleright$}%
	\begin{description}[leftmargin=15pt,itemindent=-5.2pt]}{%
	\end{description}\renewcommand{\descriptionlabel}{\descriptionlabelsave}}

\newcounter{StepsCount}				

\newcounter{StepsRefCount}

\theoremstyle{plain}
	\newtheorem{lemma}{Lemma}[section]
	\newtheorem{proposition}[lemma]{Proposition}
	\newtheorem{theorem}[lemma]{Theorem}
	\newtheorem{corollary}[lemma]{Corollary}
	\newcommand{\GenericTheoremName}{}\newtheorem{generictheorem}[lemma]{\GenericTheoremName}
\theoremstyle{definition}
	\newtheorem{definition}[lemma]{Definition}
	\newcommand{\GenericDefinitionName}{}\newtheorem{genericdefinition}[lemma]{\GenericDefinitionName}
\theoremstyle{remark}
	\newtheorem{remark}[lemma]{Remark}
	\newtheorem{example}[lemma]{Example}
	\newcommand{\GenericRemarkName}{}\newtheorem{genericremark}[lemma]{\GenericRemarkName}
\newenvironment{Theorem}{\begin{theorem}}{\par\noindent\rule{5em}{1pt}\end{theorem}}
\newenvironment{Corollary}{\begin{corollary}}{\par\noindent\rule{5em}{1pt}\end{corollary}}
\newenvironment{Definition}{\begin{definition}}{\par\noindent\rule{5em}{1pt}\end{definition}}
\newenvironment{Remark}{\begin{remark}}{\par\noindent\rule{5em}{0.5pt}\end{remark}}

\newcommand{\mc}[1]{{\mathcal{#1}}}			
\newcommand{\ms}[1]{{\mathscr{#1}}}			
\newcommand{\bb}[1]{{\mathbb{#1}}}			
\newcommand{\ov}{\overline}				
\newcommand{\To}{\longrightarrow}			

\DeclareMathOperator{\RE}{Re}				
\renewcommand{\Re}{\RE}
\DeclareMathOperator{\IM}{Im}				
\renewcommand{\Im}{\IM}

\newcommand{\Dis}[1]{${\displaystyle{#1}}$}		
\newcommand{\Side}[1]{\hfill{#1}\kern10pt}		
\newcommand{\FD}[5]{
	\DF\left\{\begin{array}{rcl}{#1}&\to &{#2}\\[#3pt] {#4}&\mapsto &{#5}\end{array}\right.}

\newcommand{\smmatrix}[4]{\Bigl(			
\begin{smallmatrix}
\hspace*{-0.2ex} #1 \hspace*{0.2ex} & \hspace*{0.2ex} #2 \hspace*{-0.2ex}
\\[0.5ex]
\hspace*{-0.2ex} #3 \hspace*{0.2ex} & \hspace*{0.2ex} #4 \hspace*{-0.2ex}
\end{smallmatrix}
\Bigr)}
\newcommand{\smfrac}[2]{{\textstyle\frac{{#1}}{{#2}}}}	

\newcommand{\Dummy}{\text{\textvisiblespace\kern1pt}}	

\DeclareMathOperator{\Id}{id}				

\newcommand{\DS}{\mid\mkern3mu}				
\newcommand{\DP}{\colon}				
\newcommand{\DF}{\colon}				
\newcommand{\DE}{\mathrel{\mathop:}=}			
\newcommand{\ED}{=\mathrel{\mathop:}}			
\newcommand{\DD}{\mkern4mu\mathrm{d}}			

\newcommand{\CL}[2]{{\ms C}({#1},{#2})}			
\newcommand{\Cl}[1]{{\ms Cl}[{#1}]}			
\newcommand{\Cln}[1]{{\ms Cl}_{\|\Dummy\|}[{#1}]}	
\newcommand{\Ca}[1]{{\ms C}_{\text{\tiny\rm\varangle}}({#1})}	
\newcommand{\CH}{\bb{CH}}				

\DeclareMathOperator{\tr}{tr}				
\newcommand{\EllOne}[1]{L^1((0,{#1}),\bb C^{2\times 2})}	
\newcommand{\EllInf}[1]{L^\infty((0,{#1}),\bb C^{2\times 2})}	
\DeclareMathOperator{\diam}{diam}			
\newcommand{\NTto}{\stackrel{\tiny\varangle}{\to}}	

\begin{document}

\begin{flushleft}
	{\Large\bf Limit behaviour of Weyl coefficients}
	\\[5mm]
	\textsc{
	Raphael Pruckner
	\,\ $\ast$\,\ 
	Harald Woracek
		\hspace*{-14pt}
		\renewcommand{\thefootnote}{\fnsymbol{footnote}}
		\setcounter{footnote}{2}
		\footnote{This work was supported by the project P\,30715--N35 of the Austrian Science Fund.
			The second author was supported by the joint project I~4600 of the Austrian
			Science Found (FWF) and the Russian Foundation of Basic Research (RFBR).}
		\renewcommand{\thefootnote}{\arabic{footnote}}
		\setcounter{footnote}{0}
	}
	\\[6mm]
	{\small
	\textbf{Abstract:}
	We study the sets of radial or nontangential limit points towards $i\infty$ of a Nevanlinna 
	function $q$. Given a nonempty, closed, and connected subset $\mc L$ of $\ov{\bb C_+}$, we explicitly construct a 
	Hamiltonian $H$ such that the radial- and outer angular cluster sets towards $i\infty$ of the Weyl coefficient $q_H$ 
	are both equal to $\mc L$. Our method is based on a study of the continuous group action of rescaling 
	operators on the set of all Hamiltonians. 
	\\[3mm]
	{\bf AMS MSC 2010:}
	34B20, 30D40, 37J99, 30J99
	\\
	{\bf Keywords:} Weyl coefficient, canonical system, cluster set, Nevanlinna function
	}
\end{flushleft}

\section{Introduction}

A Nevanlinna function is an analytic function in the open upper half-plane $\bb C_+$ 
whose values lie in $\bb C_+\cup\bb R$. Such functions are intensively studied for various reasons; 
we mention two of them.
\begin{Itemize}
\item In complex analysis they occur as regularised Cauchy-transforms of positive Poisson integrable measures, 
	e.g.\ \cite{levin:1980,kac.krein:1968a,gesztesy.tsekanovskii:2000}.
	Namely, a function $q$ is a Nevanlinna function if and only if it is of the form 
	\begin{equation}\label{C03}
		q(z)=a+bz+\int_{\bb R}\Big(\frac 1{x-z}-\frac x{1+x^2}\Big)\DD \mu(x),\quad z\in\bb C_+
		,
	\end{equation}
	where $a\in\bb R$, $b\geq 0$, and $\mu$ is a positive Borel measure on the real line with 
	$\int_{\bb R}\frac{d\mu(x)}{1+x^2}<\infty$. 
\item In spectral theory of differential operators they occur 
	as Weyl coefficients whenever H.Weyl's nested disks method is applicable, e.g.\ 
	\cite{weyl:1910,titchmarsh:1946,atkinson:1964,behrndt.hassi.snoo:2020}. 
\end{Itemize}
The connection between these two instances is that (for simplicity we suppress some technical issues and exceptional cases) 
the measure $\mu$ in the integral representation \eqref{C03} of the Weyl coefficient of an equation is a spectral measure 
for the corresponding selfadjoint model operator. 

The natural context of Weyl's method is the framework of two-dimensional canonical systems 
\begin{equation}\label{C02}
	y'(t)=zJH(t)y(t),\quad t\in(0,\infty)
	,
\end{equation}
where $z\in\bb C$ is the eigenvalue parameter, $J\DE\smmatrix 0{-1}10$, 
and the Hamiltonian $H$ of the system is assumed to satisfy $H(t)\geq 0$ and $\tr H(t)=1$ a.e., e.g.\ 
\cite{debranges:1968,hassi.snoo.winkler:2000,romanov:1408.6022v1,remling:2018}.

It is a deep theorem due to L.de~Branges that the map assigning to each Hamiltonian $H$ 
the Weyl coefficient $q_H$ of the equation \eqref{C02} is a bijection between the set of all Hamiltonians 
\begin{equation}\label{C01}
	\bb H\DE\big\{H\DF(0,\infty)\to\bb R^{2\times 2}\DS H\text{ measurable}, H(t)\geq 0, \tr H(t)=1\text{ a.e.}\big\}
\end{equation}
up to equality a.e., and the set of all Nevanlinna functions including the function identically equal to $\infty$
\[
	\mc N\DE\big\{q\DF\bb C_+\to\ov{\bb C}\DS q\text{ analytic}, q(\bb C_+)\subseteq\ov{\bb C_+}\big\}
	.
\]
Here $\ov{\bb C}$ denotes the Riemann sphere $\bb C\cup\{\infty\}$ considered as a Riemann surface in the usual way, 
and $\ov{\bb C_+}$ denotes the closure of $\bb C_+$ in the sphere, explicitly, $\ov{\bb C_+}=\bb C_+\cup\bb R\cup\{\infty\}$. 
The assignment $H\mapsto q_H$ is also called de~Branges' correspondence. 

Having available this bijection, it is a natural task to relate properties of $H$ to properties of $q_H$. 
For many properties of Hamiltonians or Nevanlinna functions it turns out to be quite involved (or even quite impossible) to find 
their counterpart on the other side of de~Branges' correspondence. 
One type of properties where some explicit relations are known is the high-energy behaviour of $q_H$, i.e., its
behaviour towards $i\infty$. It is a frequently instantiated intuition, going back at least to B.M.Levitan \cite{levitan:1952}, 
that the high-energy behaviour of $q_H$ corresponds to the local behaviour of $H$ at $0$. 
For example it is shown in \cite{eckhardt.kostenko.teschl:2018} that the nontangential limit\footnote{%
	We write $z_n\NTto i\infty$ for: $|z_n|\to\infty$ while $\arg z_n\in[\alpha,\pi-\alpha]$ for some 
	$\alpha\in(0,\frac\pi 2]$. And we write $\lim_{z\NTto i\infty}q(z)=\zeta$, if $\lim_{n\to\infty}q(z_n)=\zeta$ for 
	every sequence $z_n\NTto i\infty$. Convergence on the Riemann sphere is understood 
	w.r.t.\ the chordal metric.
	}
$\lim_{z\NTto i\infty}q_H(z)$ exists in $\ov{\bb C}$, if and only if the limit 
$\lim_{t\searrow 0}\frac 1t\int_0^t H(s)\DD s$ exists in $\bb R^{2\times 2}$. Moreover, if these limits exist, they are related 
by simple formulae. 

In this paper we investigate the situation when the Weyl coefficient does not necessarily have a limit. 
Natural substitutes for a limit value are cluster sets. We consider two variants which are fitted to nontangential approach. 
For $\alpha\in(0,\frac\pi 2]$ denote by $\Gamma_\alpha$ the Stolz angle
\[
	\Gamma_\alpha\DE\big\{z\in\bb C_+\DS\arg z\in[\alpha,\pi-\alpha]\big\}
	.
\]
\begin{Enumerate}
\item Let $M\subseteq\bb C_+$ be such that 
	\begin{equation}
	\label{C79}
		\begin{aligned}
			& M\text{ unbounded},\ \exists \alpha\in(0,\frac\pi 2]\DP M\subseteq\Gamma_\alpha,
			\\
			& \{z\in M\DS |z|\geq r\}\text{ connected for all sufficiently large }r.
		\end{aligned}
	\end{equation}
	For a Nevanlinna function $q$ we consider the cluster set 
	\[
		\CL{q}{M}\DE\big\{\zeta\in\ov{\bb C}\DS \exists z_n\in M\DP |z_n|\to\infty\wedge q(z_n)\to\zeta\big\}
		.
	\]
\item The outer angular cluster set of a Nevanlinna function $q$ is 
	\[
		\Ca{q}\DE\bigcup_{\alpha\in(0,\frac\pi 2]}\CL{q}{\Gamma_\alpha}
		=\big\{\zeta\in\ov{\bb C}\DS \exists z_n\in\bb C_+\DP z_n\NTto\infty\wedge q(z_n)\to\zeta\big\}
		.
	\]
\end{Enumerate}
We do not consider arbitrary -- possibly tangential -- approach to infinity. 

The cluster sets $\CL qM$ and $\Ca q$ are both nonempty and connected. They show different behaviour in the sense that 
$\CL qM$ is always closed, while $\Ca q$ need not have this property, cf.\ \cite{collingwood.lohwater:1966,noshiro:1960}
(see also \Cref{C96} below). 

It is known from \cite{belna.colwell.piranian:1985} (using a fractional linear transformation to pass from half-plane
to unit disk) that for every nonempty, closed, and connected subset $\mc L$ of $\ov{\bb C_+}$ there exists a Nevanlinna
function $q$ such that the radial cluster set $\CL{q}{i[1,\infty)}$ equals $\mc L$. 
In fact, in \cite[Theorem]{belna.colwell.piranian:1985} the radial boundary
interpolation problem is solved for countably many interpolation nodes, and the given solution is a Blaschke product.
Variants of this result for singular inner functions can be found in \cite[Theorem~9]{decker:1994}, or 
\cite[Theorem~3]{donaire:2001}. For smaller classes of functions, e.g.\ interpolating or thin Blaschke products, 
the radial boundary interpolation problem is in general not anymore solvable, cf.\ \cite{gorkin.mortini:2005,girela.suarez:2011}. 
The outer angular cluster set is a countable increasing union of nonempty, closed, and connected subsets, and P.M.Gauthier 
showed in \cite{gauthier:2021} (personal communication) that every set of this form can be realised as outer angular cluster set. 

Our main result in the present paper is \Cref{C70}, where we give an explicit solution to the following inverse spectral
problem:
\begin{quote}
	\emph{Given a nonempty, closed, and connected subset $\mc L$ of $\ov{\bb C_+}$, find a Hamiltonian $H$ such that 
	$\CL{q_H}M=\Ca{q_H}=\mc L$ (for arbitrary $M$ as in \cref{C79})}
\end{quote}
The Hamiltonian $H$ constructed in the proof of \Cref{C70} has the property that $q_H$ (transferred to the unit disk) 
is a Blaschke product. 

Our method of proof is based on a rescaling trick which goes back at least to Y.Kasahara \cite{kasahara:1975}, who applied it 
on the level of Krein strings, and which was exploited further in \cite{kasahara.watanabe:2010}, and in 
\cite{eckhardt.kostenko.teschl:2018} and its forthcoming extension \cite{langer.pruckner.woracek:asysupp}. 
Namely, given a Hamiltonian $H\in\bb H$, one considers rescaled Hamiltonians
\begin{equation}\label{C04}
	(\mc A_rH)(t)\DE H(\smfrac tr),\ t\in(0,\infty),\quad r>0
	.
\end{equation}
The operators $\mc A_r$ blow up the scale and thereby zoom into the vicinity of $0$. 
We will see that cluster sets of $q_H$ are related to cluster sets of the family $(\mc A_rH)_{r\geq 1}$ where the set $\bb H$ is
appropriately topologised, cf.\ \Cref{C99,C97}.
In fact, one may say that the continuous group action of rescaling operators on $\bb H$ is responsible for the 
mentioned intuition that high-energy behaviour of $q_H$ relates to local behaviour of $H$ at $0$. 

In \cite{eckhardt.kostenko.teschl:2018,kasahara:1975,kasahara.watanabe:2010} a simple continuity property of de~Branges' 
correspondence was sufficient to obtain the desired conclusions. 
This property goes back at least to \cite{debranges:1961}, where it formed a step in the existence proof of the 
inverse spectral theorem. Despite being used in the literature ever since, an explicit presentation was given only recently in 
\cite{remling:2018}. In the presently considered general situation, when limits do not necessarily exist, 
finer arguments and a thorough understanding of the topology on $\bb H$ are necessary. 

After this introduction the article is structured in three more sections. 
In Section~2 we study the appropriate topology on $\bb H$; this section is to a certain extent of expository nature.
Contrasting the presentation in \cite{remling:2018}, we introduce the topology from a higher level viewpoint. 
Namely, as an inverse limit of weak topologies on sets of Hamiltonians defined on finite intervals ($T\in(0,\infty)$)
\begin{equation}\label{C05}
	\bb H_T\DE\big\{H\DF(0,T)\to\bb R^{2\times 2}\DS H\text{ measurable}, H(t)\geq 0, \tr H(t)=1\text{ a.e.}\big\}
	.
\end{equation}
By this approach the most important features, namely compactness and metrisability, are readily built into the construction. 
Besides offering structural clarity, it also simplifies matters by avoiding the unnecessary passage 
from $L^1$ to the space of complex Borel measures made in \cite{debranges:1961,remling:2018}. 
For the convenience of the non-specialist reader, we include a complete and concise derivation of the required continuity 
of de~Branges' correspondence $H\leftrightarrow q_H$. 

In Section~3 we study the group action of rescaling operators $\{\mc A_r\DS r>0\}$ on $\bb H$, 
and relate limit points of $q_H$ with limit points of $(\mc A_rH)_{r\geq 1}$.
The case that limits exist, which has been studied in \cite{eckhardt.kostenko.teschl:2018},
is revisited in the extended preprint version of this article, 
cf.\ \cite{pruckner.woracek:limp-ASC}.

Section~4 is devoted to the proof of the main result of the paper. 
In \Cref{C70} we give the afore mentioned explicit construction of Hamiltonians whose Weyl coefficient has prescribed 
cluster set. We close the paper with stating some open problems related to \Cref{C70}.

\section{Topologising the set of Hamiltonians}

Thoroughly understanding convergence of Hamiltonians is crucial for our present investigation. 
We shall first consider Hamiltonians defined on a finite interval 
and then pass to Hamiltonians on the half-line by a limiting process. 

\subsection{Hamiltonians on a finite interval}

Recall the notation \cref{C05}:

\begin{Definition}
\label{C26}
	For $T>0$ we denote the set of all Hamiltonians on the interval $(0,T)$ by $\bb H_T$, i.e., 
	\[
		\bb H_T\DE
		\big\{H\DF(0,T)\to\bb R^{2\times 2}\DS H\text{ measurable}, H(t)\geq 0, \tr H(t)=1\text{ a.e.}\big\}
		.
	\]
	We shall always tacitly identify two Hamiltonians which coincide almost everywhere. 
\end{Definition}

\noindent
Let $\|\Dummy\|$ denote the $\ell^1$-norm on $\bb C^{2\times 2}$. 
For every positive semidefinite matrix $A=(a_{ij})_{i,j=1}^2$ it holds that $|a_{ij}|\leq\|A\|\leq 2\tr A$. 
This yields that all $H\in \bb H_T$ are entrywise (equivalently, w.r.t.\ $\|\Dummy\|$) essentially bounded by $2$. 
In particular, we have 
\[
	\bb H_T\subseteq\EllOne T
	.
\]
The space $\EllOne T$, and with it its subset $\bb H_T$, carries several natural topologies. We will work with
its norm and weak topology, $\mc T_{\|\Dummy\|_1}$ and $\mc T_w$. 

\begin{Remark}
\label{C28}
	In order to work with the weak topology, we recall the following representation of continuous functionals. 
	We have (linearly and homeomorphically)
	\begin{align*}
		\EllOne T'\cong &\, \big[L^1(0,T)^4\big]'
		\\
		\cong &\, \big[L^1(0,T)'\big]^4\cong\big[L^\infty(0,T)\big]^4\cong\EllInf T
		.
	\end{align*}
	A linear homeomorphism is given by the assignment 
	\[
		\left\{
		\begin{array}{rcl}
			\EllInf T & \to & \EllOne T'
			\\[3mm]
			(f_{ij})_{i,j=1}^2 & \mapsto & 
			\bigg[(h_{ij})_{i,j=1}^2\mapsto\sum\limits_{i,j=1}^2\int\limits_0^Th_{ij}(t)f_{ij}(t)\DD t\big)\bigg]
		\end{array}
		\right.
	\]
	Sometimes it is practical to note that $\EllOne T'$ is spanned by the set of functionals 
	\[
		\Big\{H\mapsto\int_0^Te_1^*H(t)e_2\cdot f(t)\DD t\DS 
		e_1,e_2\in\big\{{\textstyle\binom 10,\binom 01}\big\},f\in L^\infty(0,T)\Big\}
		.
	\]
\end{Remark}

\noindent
The weak topology on $\bb H_T$ has striking properties. 

\begin{lemma}
\label{C12}
	Let $T>0$. The weak topology $\mc T_w|_{\bb H_T}$ is compact and metrisable.
\end{lemma}
\begin{proof}
	Since $\bb H_T$ is uniformly bounded, it is also uniformly integrable. The Dunford-Pettis Theorem (see, e.g., 
	\cite[Theorem~4.7.18]{bogachev:2007}) yields that $\bb H_T$ is relatively compact in the weak topology of 
	$\EllOne T$. Since every $\|\Dummy\|_1$-convergent sequence has a subsequence which converges pointwise a.e., the set 
	$\bb H_T$ is $\|\Dummy\|_1$-closed. Since it is convex, it follows that it is weakly closed. Hence $\bb H_T$ is indeed
	weakly compact. 

	Since $\EllOne T$ is $\|\Dummy\|_1$-separable, the weak topology on a weakly compact subset is metrisable
	(see, e.g., \cite[Proposition~3.2.9]{fabian.etal:2001}).
\end{proof}

\noindent
We come to a variant of continuity in de~Branges' correspondence for Hamiltonians on finite intervals. 
To this end, we need some notation. 
First, denote by $\mc E$ the set of all entire $2\times 2$-matrix functions endowed with the topology $\mc T_{\rm lu}$ of 
locally uniform convergence. Second, we introduce a notation for the (transpose of the) 
fundamental solution of a canonical system.

\begin{Definition}
\label{C17}
	Let $T>0$. For $H\in\bb H_T$ we denote by $W(H;t,z)$ the unique solution of the initial value problem 
	\begin{equation}\label{C29}
		\left\{
		\begin{array}{l}
			\frac\partial{\partial t}W(H;t,z)J=zW(H;t,z)H(t),\ t\in[0,T],
			\\[2mm]
			W(H;0,z)=I,
		\end{array}
		\right.
	\end{equation}
	where $I$ is the $2\times 2$-identity matrix. 
\end{Definition}

\noindent
For every fixed $t\in[0,T]$, the matrix $W(H;t,\Dummy)$ is an entire function, i.e., $W(H;t,\Dummy)$ belongs to $\mc E$. 

\begin{Definition}
\label{C19}
	Let $T>0$. We denote by $\Psi_T$ the map 
	\[
		\Psi_T\DF\left\{
		\begin{array}{rcl}
			\bb H_T & \to & \mc E
			\\[1mm]
			H & \mapsto & W(H;T,\Dummy)
		\end{array}
		\right.
	\]
\end{Definition}

\noindent
The above announced continuity result now reads as follows. 

\begin{Theorem}[\textbf{Continuity; fundamental solution}]
\label{C20}
	\hfill\\
	Let $T>0$. Then $\Psi_T$ is $\mc T_w$-to-$\mc T_{\rm lu}$--continuous. 
\end{Theorem}

\noindent
This theorem is implicit in \cite{debranges:1961}, and, up to identification of topologies, explicit in 
\cite[Theorem~5.7]{remling:2018}. For convenience of the reader we give a complete proof. 

\begin{proof}[Proof of \Cref{C20}]
	Let 
	\begin{equation}
	\label{C18}
		W(H;t,z)=\sum_{l=0}^\infty W_l(H;t)z^l
	\end{equation}
	be the power series expansion of $W(H;t,\Dummy)$. 
	Plugging this in the equation \cref{C29}, we obtain that the coefficients $W_l(H;t)$ satisfy the reccurrance 
	\[
		W_0(H,t)=I,\quad 
		W_{l+1}(H;t)=-\int_0^t W_l(H;s)H(s)J\DD s,\ l\in\bb N
		.
	\]
	From this one inductively obtains 
	\begin{equation}
	\label{C24}
		\|W_l(H;t)\|\leq\frac{(2t)^l}{l!},\quad H\in\bb H_T,t\in[0,T],l\in\bb N
		.
	\end{equation}
	Therefore, for each compact set $K\subseteq\bb C$, the series \cref{C18} converges uniformly on 
	$\bb H_T\times[0,T]\times K$, and we have the global growth estimate 
	\[
		\|W(H;t,z)\|\leq e^{2t|z|},\quad (H,t,z)\in\bb H_T\times[0,T]\times K
		.
	\]
	Now let $(H_n)_{n\in\bb N}$ be a sequence in $\bb H_T$ which converges weakly to some $H\in\bb H_T$. 
	By uniformity in $H$ of the convergence of the series \cref{C18}, it suffices to show that 
	\[
		\forall l\in\bb N\DP \lim_{n\to\infty}W_l(H_n;T)=W_l(H;T)
	\]
	in order to conclude that $\lim_{n\to\infty}W(H_n;T,\Dummy)=W(H;T,\Dummy)$ locally uniformly on $\bb C$. 
	We use induction to show the stronger statement
	\[
		\forall l\in\bb N\DP \lim_{n\to\infty}W_l(H_n;t)=W_l(H;t) \text{ uniformly for }t\in[0,T]
		.
	\]
	For $l=0$ this is trivial. Assume that it has already been established for some $l\in\bb N$. Using the reccurrance gives 
	\begin{align}
		\|W_{l+1}(H_n;t) &\, -W_{l+1}(H;t)\|_\infty
		\nonumber
		\\
		= &\, 
		\Big\|\int_0^t W_l(H_n;s)H_n(s)J\DD s-\int_0^t W_l(H;s)H(s)J\DD s\Big\|_\infty
		\nonumber
		\\
		\leq &\, \Big\|\int_0^t\big(W_l(H_n;s)-W_l(H;s)\big)H_n(s)J\DD s\Big\|_\infty
		\nonumber
		\\
		&\, +\Big\|\underbrace{\int_0^t W_l(H;s)\big(H_n(s)-H(s)\big)J\DD s}_{\ED g_n(t)}\Big\|_\infty
		.
		\label{C21}
	\end{align}
	The first summand is estimated as 
	\[
		\Big\|\int_0^t\big(W_l(H_n;s)-W_l(H;s)\big)H_n(s)J\DD s\Big\|_\infty
		\leq T\cdot\|W_l(H_n;s)-W_l(H;s)\|_\infty\cdot 2
		,
	\]
	and tends to $0$ by the inductive hypothesis. The functions $g_n$ tend to $0$ pointwise on $[0,T]$ since 
	\[
		\|g_n(t)\|=\Big\|\int_0^T\mathds{1}_{(0,t)}(s)W_l(H;s)\cdot\big(H_n(s)-H(s)\big)\cdot J\DD s\Big\|
	\]
	and $\lim_{n\to\infty}^w H_n=H$. It holds that 
	\[
		g_n(0)=0,\quad \|g_n(t)-g_n(t')\|\leq|t-t'|\cdot\frac{(2T)^l}{l!}\cdot 4
		,
	\]
	and by the Arzela-Ascoli Theorem the family $\{g_n\DS n\in\bb N\}$ is relatively compact in $C([0,T],\bb C^{2\times 2})$. 
	Thus pointwise convergence upgrades to uniform convergence, and we obtain that also the second summand in 
	\cref{C21} tends to $0$. 
\end{proof}

\subsection{Hamiltonians on the half-line}

We turn to Hamiltonians defined on the whole half-line. Recall the notation \cref{C01}:

\begin{Definition}
\label{C33}
	We denote the set of all Hamiltonians on the half-line $(0,\infty)$ by $\bb H$, i.e., 
	\[
		\bb H\DE
		\big\{H\DF(0,\infty)\to\bb R^{2\times 2}\DS H\text{ measurable}, H(t)\geq 0, \tr H(t)=1\text{ a.e.}\big\}
		.
	\]
	Again we tacitly identify two Hamiltonians which coincide almost everywhere. 
\end{Definition}

\noindent
We consider the set of functions on the half-line as inverse limit of the sets of functions on finite intervals in the usual way. 
For $T>0$ let $\rho_T$ be the restriction map 
\[
	\rho_T\FD{\bb H}{\bb H_T}4{H}{H|_{(0,T)}}
\]
and let $\iota$ be the map 
\[
	\iota\FD{\bb H}{\prod_{T>0}\bb H_T}6{H}{(\rho_TH)_{T>0}}
\]
Then $\iota$ is injective and 
\[
	\iota(\bb H)=\Big\{(H_T)_{T>0}\in\prod_{T>0}\bb H_T\DS \forall 0<T<T'\DP H_{T'}|_{(0,T)}=H_T\Big\}
	.
\]
We use $\iota$ to pull back the topology of the product. I.e., 
we define a topology on $\bb H$ by the demand that $\iota$ becomes a homeomorphism of $\bb H$ 
onto $\iota(\bb H)$, where the codomain is topologised in the canonical way.

\begin{Definition}
\label{C83}
	Let $\mc T$ be the initial topology on $\bb H$ w.r.t.\ the one-element family $\{\iota\}$ from the product topology of
	the weak topologies on $\bb H_T$. 
\end{Definition}

\noindent
This construction automatically implies the following crucial properties. 

\begin{lemma}
\label{C82}
	The topology $\mc T$ is compact and metrisable.
\end{lemma}
\begin{proof}
	By Tychonoff's theorem and \Cref{C12} the product topology of the weak topologies on $\bb H_T$ is compact. 
	Each restriction map 
	\[
		\rho^{T'}_T\FD{\EllOne{T'}}{\EllOne T}4{F}{F|_{(0,T)}}
	\]
	is $\|\Dummy\|_1$-to-$\|\Dummy\|_1$--continuous, and hence also $w$-to-$w$--continuous. Thus $\iota(\bb H)$ is a 
	closed subset of the product, and hence also compact. 

	Consider the map 
	\[
		\kappa\FD{\prod_{T>0}\bb H_T}{\prod_{n\in\bb N}\bb H_n}6{(H_T)_{T>0}}{(H_n)_{n\in\bb N}}
	\]
	Then $\kappa$ is clearly continuous when both products are endowed with the product topology of the weak topologies. 
	Moreover, $\kappa|_{\iota(\bb H)}$ is injective. Since $\iota(\bb H)$ is compact, it is therefore a homeomorphism 
	of $\iota(\bb H)$ onto $(\kappa\circ\iota)(\bb H)$. \Cref{C12} implies that the countable product 
	$\prod_{n\in\bb N}\bb H_n$ is metrisable. It follows that $\bb H$, 
	being homeomorphic to a subspace of this product, is metrisable. 
\end{proof}

\begin{Remark}
\label{C36}
	The topology $\mc T$ constructed above coincides with the topology defined in \cite[Chapter~5.2]{remling:2018}. 
	This follows by writing out our definition and the argument which gave metrisability of $\mc T$, and remembering 
	\Cref{C28}. 

	In \cite[Proposition~2.3]{eckhardt.kostenko.teschl:2018} convergence of Hamiltonians is introduced in yet another form. 
	To see that this form coincides with convergence w.r.t.\ $\mc T$, one has to note that step functions 
	are dense in $L^1$. 
\end{Remark}

\noindent
We turn to continuity of de~Branges' correspondence. 
Recall that $\mc N$, as a subset of the space of all analytic functions of $\bb C_+$ into the Riemann sphere, 
naturally carries the topology $\mc T_{\rm lu}$ of locally uniform convergence. 

\begin{Definition}
\label{C22}
	We denote by $\Psi$ the map 
	\[
		\Psi\DF\left\{
		\begin{array}{rcl}
			\bb H & \to & \mc N
			\\
			H & \mapsto & q_H
		\end{array}
		\right.
	\]
\end{Definition}

\begin{Theorem}[\textbf{Continuity; Weyl coefficients}]
\label{C23}
	\hfill\\
	The map $\Psi$ is $\mc T$-to-$\mc T_{\rm lu}$--homeomorphic. 
\end{Theorem}

\noindent
Also this theorem is implicit in \cite{debranges:1961} and explicit in \cite[Theorem~5.7]{remling:2018}, 
and we provide a complete derivation for the convenience of the reader.

The proof of the ``finite interval variant'' \Cref{C20} relied on the uniform estimate \cref{C24} of power series 
coefficients. The proof of the present ``half-line variant'' will follow from a uniform estimate of the size of Weyl disks.

Recall that for $H\in\bb H$ and $T>0$ the Weyl disk $\Omega_{T,z}(H)$ at $z\in\bb C_+$ is the image of $\ov{\bb C_+}$ 
under the fractional linear transformation with coefficient matrix $W(H;T,z)$. Moreover, recall that the inverse stereographical 
projection is Lipschitz continuous. In fact, considering the Riemann sphere as the unit sphere whose south pole lies at the 
origin of the complex plane, the chordal distance $\chi$ of two points $\zeta,\xi\in\bb C$ 
(suppressing explicit notation of the stereographical projection) is 
\[
	\chi(\zeta,\xi)=\frac{2|\zeta-\xi|}{\sqrt{1+|\zeta|^2}\sqrt{1+|\xi|^2}}
	,
\]
and hence $\chi(\zeta,\xi)\leq2|\zeta-\xi|$, $\zeta,\xi\in\bb C\subseteq\ov{\bb C}$. 

\begin{lemma}
\label{C25}
	Let $H\in\bb H$, $T>0$, and $z\in\bb C_+$. 
	The diameter of the Weyl disk $\Omega_{T,z}$ w.r.t.\ the chordal metric can be estimated as 
	\[
		\diam_\chi\Omega_{T,z}(H)\leq\frac 8{T\cdot\Im z}
		.
	\]
\end{lemma}
\begin{proof}
	Write $H=\smmatrix{h_1}{h_3}{h_3}{h_2}$, and assume first that $\int_0^Th_2(s)\DD s\geq\frac T2$. 
	Then $\infty\notin\Omega_{T,z}(H)$. By the usual formula for the the euclidean radius of $\Omega_{T,z}(H)$, see 
	e.g. \cite[Lemma~3.11]{remling:2018}, the monotonicity result \cite[Lemma~4]{debranges:1961}, and the 
	differential equation \cref{C29}, we find 
	\[
		\diam_\chi\Omega_{T,z}(H)\leq 2\diam_{|\Dummy|}\Omega_{T,z}(H)
		\leq 2\cdot\frac 2{\Im z\cdot\int_0^T h_2(t)\DD t}\leq\frac 8{T\cdot\Im z}
		.
	\]
	Now consider the case that $\int_0^T h_2(s)\DD s<\frac T2$. 
	Then we must have $\int_0^T h_1(s)\DD s\geq\frac T2$, and the 
	already established estimate applies to $\tilde H\DE-JHJ$. 
	A computation shows that $W(\tilde H;T,z)=-JW(H;T,z)J$, and hence 
	the Weyl disk $\Omega_{T,z}(\tilde H)$ is the image of $\Omega_{T,z}(H)$ under the fractional linear transformation 
	with coefficient matrix $J$. Since $J$ is unitary, this is a rotation of the sphere, and hence isometric w.r.t.\ 
	the chordal metric. We obtain 
	\[
		\diam_\chi\Omega_{T,z}(H)=\diam_\chi\Omega_{T,z}(\tilde H)\leq\frac 8{T\cdot\Im z}
		.
	\]
\end{proof}

\begin{proof}[Proof of \Cref{C23}]
	Let ($H_n)_{n\in\bb H}$ be a sequence in $\bb H$ which converges to some $H\in\bb H$. 
	By the definition of the topology of $\bb H$, this means that $\lim_{n\to\infty}^w\rho_T(H_n)=\rho_T(H)$ for every $T>0$. 
	
	Write $W(H;T,z)=(w_{ij}(H;t,z))_{i,j=1}^2$, and denote
	\[
		Q_{n,T}(z)\DE\frac{w_{12}(H_n;T,z)}{w_{22}(H_n;T,z)},\ 
		Q_T(z)\DE\frac{w_{12}(H;T,z)}{w_{22}(H;T,z)},\quad z\in\bb C_+
		.
	\]
	Throughout the following all limits of complex numbers are understood w.r.t.\ the chordal metric $\chi$.

	Let $K\subseteq\bb C_+$ satisfy $\inf_{z\in K}\Im z>0$. 
	\Cref{C25} shows that the limit 
	\[
		q_{\tilde H}(z)=\lim_{T\to\infty}\frac{w_{12}(\tilde H;T,z)}{w_{22}(\tilde H;T,z)}
	\]
	defining the Weyl coefficient of a Hamiltonian $\tilde H$ is attained uniformly for $(\tilde H,z)\in\bb H\times K$. 
	This implies
	\begin{Itemize}
	\item $\lim_{T\to\infty}Q_{n,T}(z)=q_{H_n}(z)$ uniformly for $(n,z)\in\bb N\times K$;
	\item $\lim_{T\to\infty}Q_T(z)=q_H(z)$ uniformly for $z\in K$. 
	\end{Itemize}
	\Cref{C20} says that 
	\begin{Itemize}
	\item For each $T>0$ we have $\lim_{n\to\infty}Q_{n,T}(z)=Q_T(z)$ locally uniformly for $z\in\bb C_+$. 
	\end{Itemize}
	Together we obtain 
	\[
		q_H(z)=\lim_{T\to\infty}\lim_{n\to\infty}Q_{n,T}(z)=\lim_{n\to\infty}\lim_{T\to\infty}Q_{n,T}(z)
		=\lim_{n\to\infty}q_{H_n}(z)
	\]
	locally uniformly for $z\in\bb C_+$. 

	Being a continuous bijection of a compact space onto a Hausdorff space, $\Psi$ is a homeomorphism.
\end{proof}

\noindent
We often use continuity of $\Psi$ in another form.

\begin{Definition}
\label{C34}
	We denote by $\Phi$ the map 
	\[
		\Phi\DF\left\{
		\begin{array}{rcl}
			\bb H\times\bb C_+ & \to & \ov{\bb C_+}
			\\
			(H,w) & \mapsto & q_H(w)
		\end{array}
		\right.
	\]
\end{Definition}

\noindent
The following reformulations of continuity of $\Psi$ are obtained by elementary arguments; 
explicit proof is deferred to the preprint version \cite{pruckner.woracek:limp-ASC} of this article.

\begin{Corollary}[\textbf{Continuity; Weyl coefficients / variant}]
\label{C27}
	\hfill\\
	Each of the below properties \Enumref{1} and \Enumref{2} is equivalent to 
	$\mc T$-to-$\mc T_{\rm lu}$--continuity of $\Psi$, and hence holds true.
	\begin{Enumerate}
	\item The map $\Phi$ is continuous when $\bb H\times\bb C_+$ is endowed with the product topology of $\mc T$ 
		and the euclidean topology. 
	\item For every compact set $K\subseteq\bb C_+$ the family $\{\Phi(\Dummy,w)\DS w\in K\}$ is equicontinuous. 
	\end{Enumerate}
\end{Corollary}

\subsection{Constant Hamiltonians}

A particular role is played by Hamiltonians $H\in\bb H$ which are constant a.e.\ on $(0,\infty)$. 
We denote the set of all such as $\CH$. 

Constant Hamiltonians can be identified with the points of $\ov{\bb C_+}$. 

\begin{Definition}
\label{C15}
	Let $\Theta\DF\ov{\bb C_+}\to\CH$ be the map acting as 
	\[
		\Theta(\zeta)\DE\smmatrix{h_1}{h_3}{h_3}{h_2}
		,
	\]
	where 
	\[
		h_1\DE\frac{|\zeta|^2}{|\zeta|^2+1},\ h_2\DE\frac 1{|\zeta|^2+1},\ h_3\DE\frac{\Re\zeta}{|\zeta|^2+1}
		,
	\]
	if $\zeta\neq\infty$, and 
	\[
		\Theta(\infty)\DE\smmatrix 1000
		.
	\]
\end{Definition}

\noindent
The map $\Theta$ is bijective. Its inverse $\Theta^{-1}\DF\CH\to\ov{\bb C_+}$ is given as 
\[
	\Theta^{-1}\smmatrix{h_1}{h_3}{h_3}{h_2}=\frac{h_3+i\sqrt{h_1h_2-h_3^2}}{h_2}
	,
\]
if $h_2\neq 0$, and 
\[
	\Theta^{-1}\smmatrix 1000=\infty
	.
\]
Note that $\det\Theta(\zeta)=0$ if and only if $\zeta\in\ov{\bb R}$, and that 
$\Theta(\zeta)$ is diagonal if and only if $\zeta\in i\ov{\bb R_+}$

From the defining formulae it is obvious that for each $T>0$ the map 
$\rho_T\circ\Theta\DF\ov{\bb C_+}\to\langle\bb H_T,\mc T_{\|\Dummy\|_1}\rangle$ is continuous. Thus $\rho_T\circ\Theta$ is 
also continuous into $\mc T_w$, and hence $\Theta$ is continuous into $\langle\bb H,\mc T\rangle$. 
Since $\ov{\bb C_+}$ is compact, each of 
\[
	\langle\rho_T(\CH),\mc T_{\|\Dummy\|_1}\rangle,\quad\langle\rho_T(\CH),\mc T_w\rangle,\quad 
	\langle\CH,\mc T\rangle
\]
is homeomorphic to $\ov{\bb C_+}$. In particular, these spaces are all compact. 

\begin{Remark}
\label{C16}
	The definition of $\Theta$ is made in such a way that 
	\[
		q_{\Theta(\zeta)}(z)=\zeta,\quad z\in\ov{\bb C_+}
		,
	\]
	in other words that $\Phi(\Theta(\zeta),w)=\zeta$, $w\in\bb C_+$.
	This is shown by a simple calculation, e.g.\ \cite[\S2.2,Example~1]{eckhardt.kostenko.teschl:2018}\footnote{%
		Caution: notation in \cite{eckhardt.kostenko.teschl:2018} is different.}.
\end{Remark}

\noindent
For later use we introduce a separate notation for constant Hamiltonians corresponding to boundary points of $\ov{\bb C_+}$, 
namely, 
\[
	\CH_0\DE\Theta(\ov{\bb R})=\big\{H\in\CH\DS \det H=0\big\}
	.
\]

\section{The rescaling method}

We have already mentioned the rescaling operation $\mc A_r\DF H(\Dummy)\mapsto H(\smfrac 1r\cdot\Dummy)$ on Hamiltonians 
in \cref{C04}. In this section we put this in an appropriate framework and make a connection between cluster sets of 
$\mc A_rH$ for $r\to \infty$ and $q_H(z)$ for $z\to i\infty$. 

Clearly, $\mc A_r$ maps $\bb H$ into itself and satisfies the computation rules 
\begin{equation}
\label{C55}
	\mc A_1=\Id,\qquad \forall r,s>0\DP\mc A_r\circ\mc A_s=\mc A_s\circ\mc A_r=\mc A_{rs}
	.
\end{equation}
This just means that the map 
\begin{equation}
\label{C65}
	\left\{
	\begin{array}{rcl}
		\bb R_+\times\bb H & \to & \bb H
		\\[1mm]
		(r,H) & \mapsto & \mc A_rH
	\end{array}
	\right.
\end{equation}
is a group action of $\bb R_+$ on $\bb H$.

\begin{lemma}
\label{C60}
	The map \cref{C65} is continuous. 
\end{lemma}
\begin{proof}
	Assume we are given $H_n,H\in\bb H$ with $H_n\to H$ and $r_n,r\in\bb R_+$ with $r_n\to r$, and assume 
	without loss of generality that $\frac r2\leq r_n\leq 2r$ for all $n$. 
	We have to show that 
	\[
		\forall T>0\DP \rho_T\mc A_{r_n}H_n\stackrel{w}{\To}\rho_T\mc A_rH
		.
	\]

	Recall \Cref{C28} and let $e_1,e_2\in\big\{\binom 10,\binom 01\big\}$ and $f\in L^\infty(0,T)$ be given. 
	Denote by $\tilde f$ the extension of $f$ to $L^\infty(0,\infty)$ with $\tilde f(t)=0$, $t\geq T$. 
	Then we have 
	\begin{align*}
		& 
		\int\limits_0^T e_1^*\big((\rho_T\mc A_{r_n}H_n)(t)-(\rho_T\mc A_rH)(t)\big)e_2\cdot f(t)\DD t
		\\
		&\mkern20mu
		=\int\limits_0^T e_1^*\big(H_n(\smfrac t{r_n})-H(\smfrac t{r_n})\big)e_2\cdot f(t)\DD t
		+\int\limits_0^T e_1^*\big(H(\smfrac t{r_n})-H(\smfrac tr)\big)e_2\cdot f(t)\DD t
		\\
		&\mkern20mu
		=r_n\int\limits_0^{\frac 2r T} e_1^*\big(H_n(s)-H(s)\big)e_2\cdot\tilde f(r_ns)\DD s
		+\int\limits_0^T e_1^*\big(H(\smfrac t{r_n})-H(\smfrac tr)\big)e_2\cdot f(t)\DD t
		\\
		&\mkern20mu
		=\int\limits_0^{\frac 2r T} e_1^*\big(H_n(s)-H(s)\big)e_2\cdot\big(\tilde f(r_ns)-\tilde f(rs)\big)\DD s
		\\
		&\mkern40mu
		+\int\limits_0^{\frac 2r T} e_1^*\big(H_n(s)-H(s)\big)e_2\cdot\tilde f(rs)\DD s
		+\int\limits_0^T e_1^*\big(H(\smfrac t{r_n})-H(\smfrac tr)\big)e_2\cdot f(t)\DD t
		.
	\end{align*}
	The first summand tends to $0$ since $\|e_1^*(H_n(s)-H(s))e_2\|_\infty\leq 2$ and 
	$\|\tilde f(r_ns)-\tilde f(rs)\|_1\to 0$, the second summand since $H_n\to H$ in $\bb H$, and the third 
	since $\|H(\smfrac t{r_n})-H(\smfrac tr)\|_1\to 0$. 
\end{proof}

\noindent
The fact that \cref{C65} is a continuous group action has some immediate consequences. In our context, the following two are of
interest:

\begin{Remark}
\label{C78}
	\phantom{}
	\begin{Enumerate}
	\item For every $H\in\bb H$ and $s>0$ the map $\mc A_s$ leaves $\Cl{\mc A_rH}$ invariant. Hence, 
		\cref{C65} induces a continuous group action on the cluster set $\Cl{\mc A_rH}$. 
	\item For every $H\in\bb H$ the stabiliser 
		\[
			(\bb R_+)_H\DE\big\{r\in\bb R_+\DS \mc A_rH=H\big\}
		\]
		is a closed subgroup of $\bb R_+$. 
	\end{Enumerate}
\end{Remark}

\noindent
Item \Enumref{2} of the above remark shows that $(\bb R_+)_H$ is either equal to $\{1\}$ or $\bb R_+$, or is of the form 
$\{p^n\DS n\in\bb Z\}$ for some $p>1$. 
We have $(\bb R_+)_H=\bb R_+$ if and only if $H\in\CH$, and $(\bb R_+)_H$ is a nontrivial subgroup if and only if $H$ 
is nonconstant and multiplicatively periodic.

\begin{Remark}
\label{C84}
	The case of a nontrivial stabiliser is particularly simple:
	if $H$ is multiplicatively periodic with primitive period $p>1$, then 
	\begin{equation}
	\label{C95}
		\Cl{\mc A_rH}=\{\mc A_rH\DS 1\leq r\leq p\}
		.
	\end{equation}
	For the inclusion ``$\subseteq$'' note that $\{\mc A_rH\DS r>0\}=\{\mc A_rH\DS 1\leq r\leq p\}$, 
	and hence the orbit of $H$ 
	is compact. The reverse inclusion holds since $\mc A_{sp^n}H=\mc A_sH$ for all $n\in\bb N$ and $s>0$, and hence
	$\mc A_sH=\lim_{n\to\infty}\mc A_{sp^n}H\in\Cl{\mc A_rH}$.
\end{Remark}

\noindent
Rescaling operators have a rescaling effect on fundamental solutions. 
This is a particular case of \cite[Lemma~2.7]{eckhardt.kostenko.teschl:2018}. 
For the convenience of the reader we recall the argument. 

\begin{lemma}
\label{C58}
	Let $H\in\bb H$ and $r>0$. Then the fundamental solutions, Weyl disks, and Weyl coefficients, 
	of $H$ and $\mc A_rH$ are related as ($t\geq 0$, $z\in\bb C_+$)
	\[
		W(\mc A_r H;t,z)=W(H;\smfrac tr,rz),\ 
		\Omega_{t,z}(\mc A_rH)=\Omega_{\frac tr,rz}(H),\ 
		q_{\mc A_rH}(z)=q_H(rz)
		.
	\]
	Using the notation $\Phi$ from \Cref{C34}, the relation between Weyl coefficients writes as 
	\begin{equation}
	\label{C62}
		\forall H\in\bb H,r>0,z\in\bb C_+\DP \Phi(\mc A_rH,z)=\Phi(H,rz)
		.
	\end{equation}
\end{lemma}
\begin{proof}
	Set $\tilde W(t,z)\DE W(H;\frac tr,rz)$. Then 
	\begin{align*}
		\frac\partial{\partial t}\tilde W(t,z)J= &\, \frac 1r \frac\partial{\partial t}W(H;\frac tr,rz)
		\\
		= &\, \frac 1r\cdot rz\cdot W(H;\frac tr,rz)H\big(\frac tr\big)
		=z\tilde W(t,z)(\mc A_rH)(t)
		.
	\end{align*}
	Thus $\tilde W(t,z)$ is the fundamental solution of $\mc A_rH$. 

	The relation between Weyl disks follows immediately, and the relation between Weyl coefficients 
	follows by letting $t\to\infty$. 
\end{proof}

\noindent
The next proposition is the basis for translating cluster sets of $\mc A_rH$ to such of $q_H$. 

Given a subset $M\subseteq\bb C_+$ as in \cref{C79}, we denote the limiting directions of $M$ as 
\[
	D(M)\DE\big\{\theta\in[0,\pi]\DS \exists z_n\in M\DP |z_n|\to\infty\wedge\arg z_n\to\theta\big\}
	.
\]
Note that $D(M)$ is closed and contained in $(0,\pi)$. 

\begin{proposition}
\label{C99}
	Let $M\subseteq\bb C_+$ be as in \cref{C79} and let $H\in\bb H$. Then 
	\begin{equation}
	\label{C98}
		\CL{q_H}{M}\subseteq\Phi\big(\Cl{\mc A_rH}\times e^{iD(M)}\big)=
		\CL{q_H}{e^{iD(M)}[1,\infty)}
		.
	\end{equation}
\end{proposition}
\begin{proof}
	To show the inclusion on the left of \cref{C98}, let $w\in\CL{q_H}M$. 
	Choose $z_n\in M$ with $|z_n|\to\infty$ and $q_H(z_n)\to w$. 
	By compactness of $\bb H$ and $[0,\pi]$, we can choose a subsequence such that both limits 
	\[
		\tilde H\DE\lim_{k\to\infty}\mc A_{|z_{n_k}|}H,\quad \theta\DE\lim_{k\to\infty}\arg z_{n_k}
		,
	\]
	exist. Then $\theta\in D(M)$, and continuity of $\Phi$ implies that
	\[
		w=\lim_{k\to\infty}q_H(z_{n_k})=\lim_{k\to\infty}\Phi\big(\mc A_{|z_{n_k}|}H,e^{i\arg z_{n_k}}\big)
		=\Phi(\tilde H,e^{i\theta})
		.
	\]
	The inclusion ``$\supseteq$'' of the asserted equality on the right of \cref{C98} readily follows since $D(M)$ is 
	closed and hence 
	\[
		D\big(e^{iD(M)}[1,\infty)\big)=D(M)
		.
	\]
	To prove the reverse inclusion, let $w\in\Phi(\Cl{\mc A_rH}\times e^{iD(M)})$ be given. Write 
	$w=\Phi(\tilde H,e^{i\theta})$ with some $\tilde H\in\Cl{\mc A_rH}$ and $\theta\in D(M)$, and 
	choose $r_n\to\infty$ with $\tilde H=\lim_{n\to\infty}\mc A_{r_n}H$. Then 
	\[
		w=\Phi(\tilde H,e^{i\theta})=\lim_{n\to\infty}\Phi(\mc A_{r_n}H,e^{i\theta})
		=\lim_{n\to\infty}q_H(r_ne^{i\theta})\in\CL{q_H}{e^{i\theta[1,\infty)}}
		.
	\]
\end{proof}

\noindent
We also obtain some knowledge about outer angular cluster sets. 

\begin{proposition}
\label{C97}
	Let $H\in\bb H$. Then 
	\begin{Enumerate}
	\item \Dis{%
		\Ca{q_H}=\Phi\big(\Cl{\mc A_rH}\times\bb C_+\big)
		}
	\item \phantom{}\vspace*{-7.75mm}
		\begin{multline*}
			\mkern-18mu
			\Cl{\mc A_rH}\subseteq\CH\ \Rightarrow\ 
			\\
			\forall M\text{ as in \cref{C79}}\DP\CL{q_H}M=\Ca{q_H}=\Theta^{-1}\big(\Cl{\mc A_rH}\big)
		\end{multline*}
	\item \Dis{%
		\Cl{\mc A_rH}\cap\CH=\emptyset\ \Rightarrow\ \Ca{q_H}\text{ is open}
		}
	\end{Enumerate}
\end{proposition}
\begin{proof}
	Using \Cref{C78}\,\Enumref{1} and \cref{C98} we find
	\begin{multline*}
		\Phi\big(\Cl{\mc A_rH}\times\bb C^+\big)
		=\Phi\big(\Cl{\mc A_rH}\times e^{i(0,\pi)}\big)
		\\
		=\bigcup_{\alpha\in(0,\smfrac\pi 2]}\Phi\big(\Cl{\mc A_rH}\times e^{i[\alpha,\pi-\alpha]}\big)
		=\bigcup_{\alpha\in(0,\smfrac\pi 2]}\CL{q_H}{\Gamma_\alpha}=\Ca{q_H}
		.
	\end{multline*}
	Assume now that $\Cl{\mc A_rH}\subseteq\CH$, and set $K\DE\Theta^{-1}(\Cl{\mc A_rH})$. 
	Then 
	\[
		\Ca{q_H}=\Phi\big(\Cl{\mc A_rH}\times\bb C_+\big)
		=\mkern-5mu\bigcup_{\tilde H\in\Cl{\mc A_rH}}\mkern-20muq_{\tilde H}(\bb C_+)
		=\mkern-5mu\bigcup_{\tilde H\in\Cl{\mc A_rH}}\mkern-15mu\{\Theta^{-1}(\tilde H)\}=K
		.
	\]
	The inclusion $\CL{q_H}{M}\subseteq\Ca{q_H}$ holds trivially. 
	Let $\xi\in K$, and choose $r_n\to\infty$ with $\mc A_{r_n}H\to\Theta(\xi)$. Since $\{z\in M\DS|z|>r\}$ is connected
	for all sufficiently large $r$, we can choose for all sufficiently large $n$ points $z_n\in M$ with $|z_n|=r_n$. 
	Choose a subsequence such that $\arg z_{n_k}\to\theta$ for some $\theta\in(0,\pi)$. Then 
	\[
		q_H(z_n)=\Phi\big(\mc A_{|z_{n_k}|}H,e^{i\arg z_{n_k}}\big)\to\Phi\big(\Theta(\xi),e^{i\theta}\big)=\xi
		,
	\]
	and therefore $\xi\in\CL{q_H}{M}$. 

	Finally, assume that $\Cl{\mc A_rH}\cap\CH=\emptyset$. Then 
	\[
		\Ca{q_H}=\Phi\big(\Cl{\mc A_rH}\times\bb C_+\big)=\bigcup_{\tilde H\in\Cl{\mc A_rH}}q_{\tilde H}(\bb C_+)
		,
	\]
	and each set in the union on the right is open. 
\end{proof}

\noindent
Let us revisit the multiplicatively periodic situation.

\begin{Remark}
\label{C96}
	Let $H\in\bb H$ be nonconstant and multiplicatively periodic. Then \cref{C95} and 
	\Cref{C97}\,\Enumref{3} imply that $\Ca{q_H}$ is open. 
	In particular, the outer angular cluster set is not equal to any of the cluster sets $\CL{q_H}{M}$. 
\end{Remark}

\section{Weyl coefficients with prescribed cluster set}

In the below theorem we give an explicit construction of Hamiltonians $H$ for which the cluster set of $q_H$ 
can be computed. These Hamiltonians are piecewise constant on quickly shrinking intervals which accumulate only at the initial
point. 

In the formulation of the theorem we denote the cluster set of a sequence $(\zeta_n)_{n\in\bb N}$ in $\ov{\bb C_+}$ as 
\[
	\Cl{\zeta_n}\DE
	\big\{\zeta\in\ov{\bb C_+}\DS\exists n_k\in\bb N\DP n_k\to\infty\wedge\lim_{k\to\infty}\zeta_{n_k}=\zeta\big\}
	.
\]
Moreover, recall that $\chi$ denotes the chordal metric on $\ov{\bb C_+}$. 

\begin{Theorem}
\label{C70}
	Let $(t_n)_{n\in\bb N}$ be a sequence of positive numbers with 
	\[
		1=t_1>t_2>t_3>\ldots,\quad\lim_{n\to\infty}t_n=0,\quad\lim_{n\to\infty}\frac{t_{n+1}}{t_n}=0
		,
	\]
	and let $(\zeta_n)_{n\in\bb N}$ be a sequence of points on $\ov{\bb C_+}$ with 
	\[
		\lim_{n\to\infty}\chi(\zeta_{n+1},\zeta_n)=0
		.
	\]
	Define $H$ to be the piecewise constant Hamiltonian 
	\begin{equation}
	\label{C80}
		H(t)\DE
		\begin{cases}
			\Theta(\zeta_n) &\hspace*{-3mm},\quad t\in(t_{n+1},t_n],n\in\bb N,
			\\[1mm]
			\Theta(0) &\hspace*{-3mm},\quad t\in(1,\infty).
		\end{cases}
	\end{equation}
	Then, for every $M$ as in \cref{C79},
	\[
		\Ca{q_H}=\CL{q_H}{M}=\Cl{\zeta_n}
		.
	\]
\end{Theorem}

\noindent
An elementary argument shows that for every nonempty, closed, and connected subset $\mc L$ of $\ov{\bb C_+}$ there exists a
sequence $(\zeta_n)_{n\in\bb N}$ with $\Cl{\zeta_n}=\mc L$ (an explicit proof can be found in \cite{pruckner.woracek:limp-ASC}).
Thus we obtain an explicit solution of an inverse problem dealing with boundary interpolation. 

\begin{Corollary}
\label{C85}
	Let $\mc L\subseteq\ov{\bb C_+}$ be nonempty, closed, and connected. Then we can construct a Hamiltonian 
	$H$ for whose Weyl coefficient $q_H$ the outer angular and radial cluster sets at $i\infty$ are both equal to $\mc L$. 
\end{Corollary}

\noindent
We turn to the proof of \Cref{C70}. The crucial step is presented in the next lemma. 
Here we denote by $d_{\|\Dummy\|}$ the metric induced by the $L^1$-norm.

\begin{lemma}
\label{C93}
	Let $H\in\bb H$ and assume that 
	\begin{equation}
	\label{C94}
		\lim_{r\to\infty} d_{\|\Dummy\|}\big(\rho_1\mc A_rH,\rho_1\Theta(\CH)\big)=0
		.
	\end{equation}
	Moreover, denote 
	\begin{equation}
	\label{C86}
		K\DE\big\{\xi\in\ov{\bb C_+}\DS \rho_1\Theta(\xi)\in\Cln{\rho_1\mc A_rH}\big\}
		.
	\end{equation}
	Then 
	\begin{Enumerate}
	\item \Dis{%
		\forall T>0\DP \Cln{\rho_T\mc A_rH}=\rho_T\Theta(K)
		,
		}
	\item \Dis{%
		\Cl{\mc A_rH}=\Theta(K)
		.
		}
	\end{Enumerate}
\end{lemma}
\begin{proof}
	Let $T>0$ and set 
	\[
		K_T\DE\big\{\xi\in\ov{\bb C_+}\DS \rho_T\Theta(\xi)\in\Cln{\rho_T\mc A_rH}\big\}
		.
	\]
	The relation 
	\begin{align*}
		\|\rho_T\mc A_rH-\rho_T\Theta(\xi)\|_1= &\, \int_0^T\|H(\smfrac tr)-\Theta(\xi)\|\DD t
		\\
		= &\, T\int_0^1\|H(\smfrac Tr\cdot t)-\Theta(\xi)\|\DD t=T\|\rho_1\mc A_{\frac rT}H-\rho_1\Theta(\xi)\|_1
		,
	\end{align*}
	which holds for all $H\in\bb H$ and $\xi\in\ov{\bb C_+}$, shows that 
	\begin{align}
		\label{C92}
		\liminf_{r\to\infty}\|\rho_T\mc A_rH-\rho_T\Theta(\xi)\|_1
		= &\, 
		T\cdot\liminf_{r\to\infty}\|\rho_1\mc A_{\frac rT}H-\rho_1\Theta(\xi)\|_1
		,
		\\
		\label{C91}
		d_{\|\Dummy\|}\big(\rho_T\mc A_rH,\rho_T\Theta(\CH)\big)
		= &\, 
		T\cdot d_{\|\Dummy\|}\big(\rho_1\mc A_{\frac rT}H,\rho_1\Theta(\CH)\big)
		.
	\end{align}
	The relation \cref{C92} implies that $K_T=K$ for all $T>0$, and \cref{C91} that 
	\begin{equation}
	\label{C90}
		\forall T>0\DP
		\lim_{r\to\infty} d_{\|\Dummy\|}\big(\rho_T\mc A_rH,\rho_T\Theta(\CH)\big)=0
		.
	\end{equation}
	Since $\rho_T\Theta(\CH)$ is compact w.r.t.\ $\|\Dummy\|_1$, and hence closed, \cref{C90} in turn implies that 
	\[
		\Cln{\rho_T\mc A_rH}\subseteq\rho_T\Theta(\CH)
		.
	\]
	Item \Enumref{1} of the present assertion follows. 

	The inclusion ``$\supseteq$'' in item \Enumref{2} holds because of a general argument. Namely, it holds 
	for every Hamiltonian $H\in\bb H$ that 
	\[
		\big\{\tilde H\in\bb H\DS\forall T>0\DP\rho_T\tilde H\in\Cln{\rho_T\mc A_rH}\big\}
		\subseteq\Cl{\mc A_rH}
		.
	\]
	To show this, assume that $\tilde H$ belongs to the set on the left. We choose inductively numbers $r_n>0$, such that 
	\[
		\forall n\in\bb N\DP r_{n+1}\geq r_n+1\wedge\|\rho_n\mc A_{r_n}H-\rho_n\tilde H\|_1\leq\frac 1n
		.
	\]
	Given $T>0$, we have for all $n\geq T$ 
	\[
		\|\rho_T\mc A_{r_n}H-\rho_T\tilde H\|_1\leq\|\rho_n\mc A_{r_n}H-\rho_n\tilde H\|_1\leq\frac 1n
		,
	\]
	and hence $\rho_T\mc A_{r_n}H\stackrel{\|\Dummy\|_1}{\To}\tilde H$. This clearly implies that $\mc A_{r_n}H\to\tilde H$. 

	The reverse inclusion ``$\subseteq$'' in item \Enumref{2} relies on the assumption \cref{C94}. Assume that 
	$\tilde H\in\Cl{\mc A_rH}$ and choose a sequence $r_n\to\infty$ such that $\mc A_{r_n}H\to\tilde H$. Then 
	\begin{equation}
	\label{C89}
		\forall T>0\DP \rho_T\mc A_{r_n}H\stackrel{w}{\To}\tilde H
		.
	\end{equation}
	Let $T>0$. Since \cref{C90} holds and $\rho_T\Theta(\CH)$ is compact w.r.t.\ $\|\Dummy\|_1$, we find a point 
	$\xi\in\ov{\bb C_+}$ and a subsequence $(r_{n_k})_{k\in\bb N}$ (both depending on $T$), such that 
	\begin{equation}
	\label{C88}
		\rho_T\mc A_{r_{n_k}}H\stackrel{\|\Dummy\|_1}{\To}\Theta(\xi)
		.
	\end{equation}
	Together with \cref{C89} we see that $\rho_T\tilde H=\rho_T\Theta(\xi)$. 
	It follows that $\xi$ is independent of $T$ and that $\tilde H=\Theta(\xi)$. By \cref{C88}, used for $T=1$, we have 
	$\xi\in K$. 
\end{proof}

\begin{proof}[Proof of \Cref{C70}]
	The function 
	\[
		\rho_1\circ\Theta\DF \ov{\bb C_+}\to L^1((0,1),\bb C^{2\times 2})
	\]
	is $\chi$-to-$\|\Dummy\|_1$--continuous and injective. Since $\ov{\bb C_+}$ is compact, it is therefore uniformly
	continuous and a homeomorphism onto its image. Let $\omega\DF\bb R_+\to\bb R_+$ be the modulus of continuity of 
	$\rho_1\circ\Theta$, so that 
	\[
		\lim_{\delta\to 0}\omega(\delta)=0\ \wedge\ 
		\forall\zeta,\xi\in\ov{\bb C_+}\DP \|\rho_1\Theta(\zeta)-\rho_1\Theta(\xi)\|_1\leq\omega(\chi(\zeta,\xi))
		.
	\]
	We show that 
	\begin{equation}
	\label{C87}
		\forall n\in\bb N,r\in\big[\smfrac 1{t_n},\smfrac 1{t_{n+1}}\big]\DP
		\|\rho_1\mc A_rH-\rho_1\Theta(\zeta_n)\|_1\leq 4\frac{t_{n+2}}{t_{n+1}}+\omega(\chi(\zeta_{n+1},\zeta_n))
		.
	\end{equation}
	To see this, estimate
	\begin{align*}
		&
		\int_0^{rt_{n+2}}\|H(\smfrac tr)-\Theta(\zeta_n)\|\DD t\leq 4rt_{n+2}\leq 4\frac{rt_{n+2}}{rt_{n+1}}
		,
		\\
		&
		\int_{rt_{n+2}}^{rt_{n+1}}\|H(\smfrac tr)-\Theta(\zeta_n)\|\DD t
		=\int_{rt_{n+2}}^{rt_{n+1}}\|\Theta(\zeta_{n+1})-\Theta(\zeta_n)\|\DD t
		\\
		&\mkern210mu
		\leq \|\rho_1\Theta(\zeta_{n+1})-\rho_1\Theta(\zeta_n)\|_1\leq\omega(\chi(\zeta_{n+1},\zeta_n))
		,
		\\
		&
		\int_{rt_{n+1}}^1\|H(\smfrac tr)-\Theta(\zeta_n)\|\DD t=0
		.
	\end{align*}
	For $r\geq 1$ let $n(r)\in\bb N$ be the unique number with $[\smfrac 1{t_n},\smfrac 1{t_{n+1}})$. Then 
	$\lim_{r\to\infty}n(r)=\infty$. The right side of \cref{C87} tends to $0$ when $n$ tends to $\infty$, and hence 
	for every sequence $r_k\to\infty$ we have 
	\[
		\lim_{k\to\infty}\|\rho_1\mc A_{r_k}H-\rho_1\Theta(\zeta_{n(r_k)})\|_1=0
		.
	\]
	This shows that \cref{C94} holds and that 
	\[
		\Cln{\rho_1\mc A_rH}\subseteq\Cln{\rho_1\Theta(\zeta_n)}
		.
	\]
	If $n_k\to\infty$, then \cref{C87} shows that 
	\[
		\lim_{k\to\infty}\|\rho_1\mc A_{\frac 1{t_{n_k}}}H-\rho_1\Theta(\zeta_{n_k})\|_1=0
		,
	\]
	and it follows that 
	\[
		\rho_1\Theta\big(\Cl{\zeta_n}\big)\subseteq\Cln{\rho_1\mc A_rH}
		.
	\]
	Since $\rho_1\circ\Theta$ is a homeomorphism between compact sets, 
	\[
		\rho_1\Theta\big(\Cl{\zeta_n}\big)=\Cln{\rho_1\Theta(\zeta_n)}
		,
	\]
	and we see that the set $K$ from \cref{C86} is equal to $\Cl{\zeta_n}$. 

	The asserted properties of $q_H$ now follow from \Cref{C93} and \Cref{C97}\,\Enumref{2}.
\end{proof}

\noindent
Let us pass from half-plane to unit disk with the fractional linear transformation $\beta(z)\DE\frac{z-i}{z+i}$, 
which maps $\ov{\bb C^+}$ onto the closed unit disk $\ov{\bb D}$ with $\beta(\infty)=1$. 

\begin{Remark}
\label{C81}
	Consider a Hamiltonian of the form \cref{C80}. Since $H$ is constant equal to $\Theta(0)$ on the interval $(1,\infty)$, 
	the Weyl coefficient $q_H$ is given as $q_H=\frac{w_{12}}{w_{22}}$, where 
	\[
		w_{12}(z)\DE(1,0)W(H;1,z)\binom 01,\quad w_{22}(z)\DE(0,1)W(H;1,z)\binom 01
		.
	\]
	Since $\det H$ is constant equal to $0$, the entire function $W(H;1,z)$ is of zero exponential type. 
	The function 
	\[
		B(z)\DE\beta\circ q_H\circ\beta^{-1}=\frac{w_{12}-iw_{22}}{w_{12}+iw_{22}}\circ\beta^{-1}
	\]
	is thus a Blaschke product whose zeroes have no finite accumulation point. 

	Cluster sets of $q_H$ towards $i\infty$ clearly correspond to cluster sets of $B$ towards $1$. 
	Thus we reobtain the fact that for every nonempty, closed, and connected subset $\mc L$ of $\ov{\bb D}$, 
	there exists a Blaschke product whose outer angular and radial cluster sets at $1$ are equal to $\mc L$.
\end{Remark}

\noindent
In the context of the present construction and its consequences for functions on the disk 
some open questions occur:
\begin{Enumerate}
\item We do not know if the function constructed in the above way has also cluster set $\mc L$ when $z$ is 
	allowed to approach $1$ in an unrestricted, possibly tangential, way. 
\item We do not know if our construction method can be modifed so to obtain results about simultaneous boundary 
	interpolation at more than one point (as done for radial cluster sets in 
	\cite{belna.colwell.piranian:1985,decker:1994,donaire:2001}).
\item We do not know if our construction method can be modified to produce approach to cluster values 
	along a prescribed curve when $z$ approaches the point $1$ radially (as in \cite[Theorem~1]{donaire:2001}).
\item We do not know an analogue of \Cref{C70} for outer angular cluster sets which realises any countable increasing union of
	nonempty closed connected sets (as in \cite{gauthier:2021}).
\end{Enumerate}
Concerning the third question we have some preliminary results indicating that the answer is affirmative.

{\footnotesize
\begin{flushleft}
	R.\,Pruckner\\
	Institute for Analysis and Scientific Computing\\
	Vienna University of Technology\\
	Wiedner Hauptstra{\ss}e\ 8--10/101\\
	1040 Wien\\
	AUSTRIA\\
	email: raphael.pruckner@tuwien.ac.at\\[5mm]
\end{flushleft}
}
{\footnotesize
\begin{flushleft}
	H.\,Woracek\\
	Institute for Analysis and Scientific Computing\\
	Vienna University of Technology\\
	Wiedner Hauptstra{\ss}e\ 8--10/101\\
	1040 Wien\\
	AUSTRIA\\
	email: \texttt{harald.woracek@tuwien.ac.at}\\[5mm]
\end{flushleft}
}

\end{document}